\documentclass{article}
\usepackage{amsfonts,amsmath,amsthm,amssymb}
\usepackage{graphicx, color}
\usepackage[margin=1.2in]{geometry}
\usepackage[colorlinks=true,urlcolor=blue]{hyperref} 

\newtheorem{theorem}{Theorem}[section]
\newtheorem{lemma}[theorem]{Lemma}
\newtheorem{proposition}[theorem]{Proposition}
\newtheorem{corollary}[theorem]{Corollary}
\newtheorem{definition}[theorem]{Definition}

\newcommand{\refT}[1]{Theorem~\ref{#1}}
\newcommand{\refC}[1]{Corollary~\ref{#1}}
\newcommand{\refL}[1]{Lemma~\ref{#1}}
\newcommand{\refD}[1]{Definition~\ref{#1}}
\newcommand{\refS}[1]{Section~\ref{#1}}
\newcommand{\refP}[1]{Proposition~\ref{#1}}
\newcommand{\refF}[1]{Figure~\ref{#1}}
\newcommand{\refand}[2]{\ref{#1} and~\ref{#2}}

\newcommand\cA{\mathcal A}
\newcommand\cB{\mathcal B}
\newcommand\cD{\mathcal D}
\newcommand\cF{\mathcal F}
\newcommand\cH{\mathcal H}
\newcommand\cI{\mathcal I}
\newcommand\cM{\mathcal M}
\newcommand\cN{\mathcal N}
\newcommand\cP{\mathcal P}
\newcommand\cR{{\mathcal R}}
\newcommand\cS{{\mathcal S}}

\newcommand{\E}[1]{{\mathbb E}\left[#1\right]}
\newcommand{\p}[1]{{\mathbb P}\left(#1\right)}
\newcommand{\I}[1]{{\mathbf 1}_{[#1]}}

\newcommand{\Bin}[1]{\mathrm{Bin}(#1)}
\newcommand{\Poi}[1]{\mathrm{Poi}(#1)}

\newcommand{\N}{\mathbb{N}}
\newcommand{\R}{\mathbb{R}}
\newcommand{\Z}{\mathbb{Z}}
\newcommand{\eps}{\varepsilon}
\newcommand{\todist}{\stackrel{\mathcal{L}}{\longrightarrow}}
\newcommand{\eqdist}{\stackrel{\mathcal{L}}{=}}
\newcommand{\floor}[1]{\lfloor #1 \rfloor}
\newcommand{\n}{{(n)}}
\newcommand{\nl}{{(n_\ell)}}

\title{Depth of vertices with high degree in random recursive trees}

\author{Laura Eslava\thanks{Universidad Nacional Autonoma Mexico, Instituto de investigaciones en matematicas aplicadas y en sistemas, Mexico; e-mail: {\tt laura@sigma.iimas.unam.mx}. }}

\begin{document}

\maketitle





\begin{abstract}
  Let $T_n$ be a random recursive tree with $n$ nodes. List vertices of $T_n$ according to a decreasing order of their degrees as $v^{(1)},\ldots,v^{(n)}$, and write $\deg(v^{(i)})$ and $h(v^{(i)})$ for the degree of $v^{(i)}$ and the distance of $v^{(i)}$ from the root, respectively. We prove that, as $n \to \infty$ along suitable subsequences, 
\[
\bigg(\deg(v^{(i)}) - \lfloor \log_2 n \rfloor, \frac{h(v^{(i)}) - \mu\log_e n}{\sqrt{\sigma^2\log_e n}}\bigg) \to ((P_i,i \ge 1),(N_i,i \ge 1))\, ,
\]
where $\mu=1-(\log_2 e)/2$, $\sigma^2=1-(\log_2 e)/4$, $(P_i,i \ge 1)$ is a Poisson point process on $\mathbb{Z}$ and $(N_i,i \ge 1)$ is a sequence of independent standard Gaussians. We additionally establish joint normality for the depths of finitely many uniformly random vertices in $T_n$, which extends results from Devroye in 1988 and Mahmoud in 1991. The joint limit holds even if the random vertices are conditioned to have large degree; in particular, both the mean and variance of the conditional depths remain of orden $\ln n$.

Our results are based on a $n!$-to-1 correspondence between (a representation of) Kingman's $n$-coalescent and random recursive trees; a utility that was observed in work by Pittel 1994 and recovered by Addario-Berry and the author in 2018. 
\end{abstract}

\maketitle

\section{Introduction}

Random recursive trees (RRTs) have been widely studied since their introduction, \cite{NaRapoport70}, and are closely related to binary search trees, preferential attachment trees and increasing trees in general, see e.g. \cite{BergeronFlajoletSalvy92,Drmota}.  For $n\ge 1$, let $[n]:=\{1,\ldots, n\}$ and let $\cI_n$ be the class of increasing trees  with vertex set $[n]$; that is, rooted, labelled trees where labels are increasing along paths starting at the root. 
A random recursive tree on $n$ vertices is uniformly distributed on $\cI_n$ and may be regarded as a degenerate case of linear preferential attachment trees. We will be concerned with the interplay  between  degree  and  depth  of  vertices in RRTs and will contrast our results with known results on linear preferential attachment trees \cite{Bhamidi07,Mori05}.

For $\alpha\in [0,\infty)$, the linear preferential attachment process $(T_{\alpha, n};\, n\ge 1)$ is defined as follows. Let $T_{\alpha,1}$ be a single vertex labeled 1. For $n\in \N$ let $T_{\alpha,n+1}$ be the tree obtained from $T_{\alpha,n}$ by adding an edge from a new vertex labelled $n+1$ to vertex $v_n\in [n]$; where, conditional on $T_{\alpha,n}$, $\p{v_n=v}$ is proportional to $\alpha\deg_{T_{\alpha, n}}(v)+1$. Let
\begin{align}
\Delta_{\alpha, n}:=& \max_{v\in [n]} \deg_{T_{\alpha,n}}(v), \label{dfn:Dn}\\
\cM_{\alpha,n}:=&    \{v\in [n]:\, \deg_{T_{\alpha,n}}=\Delta_{\alpha,n}\}. \label{dfn:Mn}
\end{align}

Random recursive trees correspond to $\alpha=0$ and so the choice of each $v_n$ is independent from the past and uniformly random on $[n-1]$. In what follows, we omit the index $\alpha$ from the notation of RRTs. 

The maximum degree of a RRT satisfies $\Delta_n /\log n\to 1$ a.s. as $n\to \infty$ \cite{DevroyeLu95} and the limiting law of $\Delta_n-\floor{\log n}$ has, up to lattice effects, a (discrete) Gumbel distribution. 
This was first shown using singularity analysis of generating functions in \cite{GohSchmutz02}. \cite{AddarioEslava17} provides a probabilistic proof, while \cite{Eslava21} uses the Chen-Stein method to improve the tail bounds on the limiting distribution. 
For $\alpha>0$, it has been proven that $\Delta_{\alpha, n}/n^{1/(2+1/\alpha)}$, converges a.s. and in $L_p$ to a positive, finite random variable with absolutely continuous distribution, \cite{Mori05}. Furthermore, almost surely there is a unique vertex attaining the maximum degree and its label $v^*$ is finite;  that is, $\cM_{\alpha,n}=\{v^*\}$ and $v^*=O_p(1)$ a.s. \cite{Bhamidi07}. Recently, the case of maximal degree vertices in a wider class of growing networks is analysed in \cite{BanerjeeBhamidi20+}.

In this work we are concerned with the depths of vertices with extremal degree values. We say that $v$ is a \emph{high-degree vertex} if $\deg_{T_n}(v)\ge c\ln n$ for some $c>0$ and that $v$ is a \emph{near-maximum degree vertex} if $|\deg_{T_n}(v)-\log n|<o(\ln n)$.
Addario-Berry and the author describe the number of high-degree vertices in $T_n$ via the sequence $(\deg_{T_n}(v^{(i)})-\floor{\log n},\, v^{(i)}\in [n])$, where $(v^{(i)})_{i\in [n]}$ are listed in decreasing order of their degrees; see \cite[Theorem 1.2]{AddarioEslava17}. They show that, along suitable subsequences, this sequence converges in distribution to a Poisson point process $\cN$ in $\Z\cup \{\infty\}$ with $\E{|\cN\cap [j,\infty]|}=\Theta(2^{-j})$ for all $j\in \Z$.

The main result of this paper extends such convergence to consider not only the degrees of the vertices but also their depths; see \refT{thm:joint conv}. A remarkable consequence of \refT{thm:joint conv}, in contrast to the behaviour of preferential attachement trees, is that for RRTs there is a random number of maximum-degree vertices and their labels are increasing with $n$; see \refC{cor:max conv}. On the other hand, \refT{thm:joint conv} is essentially based on \refT{thm:depths conv} which states that the depths of any finite collection of uniformly chosen vertices, conditional on exceeding any given degree, are asymptotically normal and independent; this extends the known limiting distribution of the depth of a uniformly chosen vertex in $T_n$ \cite{Devroye88, Mahmoud91}. In particular, the distributional convergence holds when conditioning the vertices on being high-degree or near-maximum-degree vertices. 

The central technique of this paper is the  link between RRTs and a representation of Kingman's coalescent developed in \cite{ab14, AddarioEslava17,Eslava21} and it is based on the fact that bijections preserve the uniform measure on finite probability spaces. In this case, there is an $n!$-to-1 mapping between the space on $n$-coalescents (described as $n$-chains in \refS{sec:Kingman}) and increasing trees on $n$ vertices. For a general description of Kingman's coalescent, see \cite[Chapter~2]{Berestycki09} and for a comparison between distinct representations of Kingman's coalescent, see \cite{Eslava21}. 

In contrast to standard constructions, our representation of Kingman's coalescent is based on the data structure known as union-find trees. In fact, the (standard) binary tree representation of Kingman's coalescent had been previously used to study union-find trees, \cite{Devroye88}, and the connection between the results therein and the height of RRTs is mentioned in \cite{Pittel94}. Although the connection between Kingman's coalescent and random recursive trees had been observed, prior to the works in \cite{AddarioEslava17,Eslava21}, its utility in studying vertex degrees seems to have gone unremarked. 

We underline two aspects of the correspondence between Kingman's coalescent and RRTs. First, it allows us to shift the perspective of vertices arriving at distinct times (and their depths and degrees being dependent on such arrival times) to a perspective where all vertices are exchangeable. Second, informally described, we may associate a random number of coin flips to each vertex in the coalescent, then its degree corresponds to the length of the first streak of heads while its depth  corresponds to the total number of tails (see  \refS{sec:Kingman} for more details). The crux of the analysis for our results then relies on controlling the correlations that arise in considering several vertices at once. 

{\bf Notation.}
We denote natural logarithms and logarithms with base 2 with $\ln(\cdot)$ and $\log (\cdot)$, respectively. We let $\N:=\{1,2,\ldots\}$ denote the natural number and use $\Z_{\ge 0}:=\{0,1,\ldots\}$ to include zero. The cumulative distribution function of the standard Gaussian distribution is denoted $\Phi$; when there is no ambiguity, we also use $\Phi$ for the measure of the distribution. We write $\cB_\R$ for the Borel sets of $\R$. For $x\in \R$ and $k\in \N$ let $(x)_k=x(x-1)\cdots (x-k+1)$ and write $(x)_0=1$. 

We consider labeled rooted trees $T=(V,E)$ where edges are naturally directed towards the root which is denoted $r(T)$. That is, we write $e=uv$ for an edge with tail $u$ and head $v$. For a vertex $v\in V(T)$, we write $\deg_T(v)$ for the number of edges directed towards $v$ in $T$ and let $h_T(v)$ denote the distance between $v$ and $r(T)$; we call $\deg_T(v)$ and $h_T(v)$ the degree and depth of $v$ respectively. Finally, a forest $F$ is a set of trees whose vertex sets are pairwise disjoint. Denote by $V(F)$ and $E(F)$, respectively, the union of the vertex and edge sets of the trees contained in $F$ and for a vertex $v\in V(F)$, we let $\deg_F(v)$ and $h_F(v)$ be the degree and depth of $v$ in the tree in $F$ that contains $v$. 

\subsection{Statement of results}

Our first result establishes the asymptotic independence of the depth of $k$ distinct vertices uniformly chosen at random. The distribution of their depths is asymptotically normal and independent, even after conditioning on arbitrary lower bounds for the vertices' degrees.

\begin{theorem}\label{thm:depths conv}
Let $k\in \N$. Let $(N_i)_{i\in [k]}$ be i.i.d.~standard Gaussian variables and let $(v_i)_{i\in [k]}$ be $k$ distinct vertices in $T_n$ chosen uniformly at random. For every $(a_1,\ldots,a_k)\in [0,2)^k$ and $(b_1,\ldots, b_k)\in \Z^k$, the conditional law of 
\[ \left( \frac{h_{T_n}(v_i)-(1-a_i/2) \ln n}{\sqrt{(1-a_i/4) \ln n}},\, i\in [k] \right)_,\]
given that $\deg_{T_n} (v_i) \ge \floor{a_i\ln n} +b_i$ for all $i\in [k]$, converges to the law of $(N_i)_{i\in [k]}$.
\end{theorem} 

The case of \refT{thm:depths conv} with $k=1$ and $a_1=b_1=0$ is implicitly established in \cite{Devroye88, Mahmoud91} where, $h_{T_n}(n)$, the insertion depth in a RRT is studied. 

For the next results, we assume that $(n_\ell)_{\ell\ge 1}$ satisfies $\log n_\ell -\floor{\log n_\ell}\to \eps$  as $\ell \to \infty$, for some $\eps\in [0,1]$.
This condition is necessary due to a lattice effect caused by the fact that $(\deg_{T_n}(v))_{v\in [n]}$ are integer valued. 
Let $\cP$ be a Poisson point process on $\R$ with intensity $\lambda(x)=2^{-x}\ln 2$. Since $\cP\cap [0,\infty)\eqdist \Poi{1}$ there exists a well defined ordering $\cP=\{P_1,P_2,\ldots\}$ such that 
\begin{align}\label{dfn:P}
|\cP\cap [P_i,\infty)|=|\cP\cap (P_i,\infty)|+1=i.    
\end{align}
For $n\ge 1$, list vertices of $T_n$ in decreasing order of degree as $v^{(1)}, \ldots, v^{(n)}$ (breaking ties uniformly at random). 

\begin{theorem}\label{thm:joint conv}
Let $(N_i)_{i\ge 1}$ be i.i.d.~standard Gaussian variables, independent of $\cP$ as defined above. For $\eps\in [0,1]$ and any increasing sequence of integers $(n_\ell)_{\ell\ge 1}$ for which $\log n_\ell -\floor{\log n_\ell}\to \eps$ as $\ell \to \infty$,  then the sequence
\[
\left(\left(\deg_{T_{n_\ell}}(v^{(i)}) - \floor{ \log n_\ell }, \frac{h_{T_{n_\ell}}(v^{(i)}) - (1-(\log e)/2) \ln n_\ell}{\sqrt{(1-(\log e)/4)\ln n_\ell}}\right), {i\in[n]}\right) \]
converges in law to $((\floor{P_i+\eps},N_i),i \ge 1))$.
\end{theorem}

For the next corollary of \refT{thm:joint conv}, recall the definition of $\cM_{n_\ell}$ in \eqref{dfn:Mn}. For $\eps>0$ let $\mathbf M_\eps$ be defined, for each $k\ge 1$, by
\begin{equation}\label{dfn:Me}
\p{\mathbf M_{\eps}=k}:=\sum_{m\in \Z} e^{-2^{-m+\eps}} \frac{2^{-(m+1-\eps)k}}{k!};
\end{equation}
then $\mathbf M_\eps$ is a mixture of Poisson random variables conditioned to being positive; see \refL{lem:dist M}. The following result focuses on
vertices attaining the maximum degree. 
 
\begin{corollary}\label{cor:max conv}
Let $(N_i)_{i\ge 1}$ be i.i.d.~standard Gaussian variables. For  $\eps\in [0,1]$ and any increasing sequence of integers $(n_\ell)_{\ell\ge 1}$ for which $\log n_\ell-\floor{\log n_\ell}\to \eps$ as $\ell\to \infty$, then 
\[\left(\frac{h_{T_{n_\ell}}(i)-(1-(\log e)/2) \ln n_\ell}{\sqrt{(1-(\log e)/4) \ln n_\ell}},\, i\in \cM_{n_\ell} \right) \todist (N_i,\, 1\le i\le \mathbf M_\eps).\]
\end{corollary}

In other words, the depths of vertices attaining the maximum degree are asymptotically normal and independent. Furthermore, the number $|\cM_n|$ of such vertices converges to a mixture of Poisson random variables conditioned to be strictly positive.

\subsection{Paper overview} 

In the next section we provide an overview of the proof of the main theorems. \refT{thm:depths Kingman} below is an equivalent statement to \refT{thm:depths conv} in terms of Kingman's coalescent whose construction is given in \refS{sec:Kingman}; in addition, we lay out the main intermediate results (\refP{prop:moments} and Lemmas \ref{lem:skorohod}--\ref{lem:neglig depth}) for the proof of \refT{thm:depths Kingman}. \refS{sec:main} is then concerned with proving Lemmas \refand{lem:cond depth}{lem:neglig depth} while \refS{sec:Proof moments} contains the proof of \refP{prop:moments}. The proofs of \refT{thm:joint conv} and \refC{cor:max conv} are given in \refS{sec:Main proofs}.

\section{Outline for the proof of Theorems \refand{thm:depths conv}{thm:joint conv}}\label{sec:Outline}


Let $T_n=([n],E_n)$ be a RRT and let $\sigma:[n]\to [n]$ be an independent, uniformly random permutation. Define $T'_n=([n],E'_n)$ where 
\begin{align}\label{dfn:En}
E'_n=\{\sigma(u)\sigma(v):\, uv\in E_n\};   
\end{align}
that is, $T'_n$ is obtained from $T_n$ by relabelling its vertices with a uniformly random permutation. 
By definition of $T'_n$, \refT{thm:depths Kingman} below is equivalent to \refT{thm:depths conv}; its proof overview is detailed in \refS{sec:Kingman}.

\begin{theorem}\label{thm:depths Kingman}
Let $k\in \N$. Let $(N_i)_{i\in [k]}$ be i.i.d.~standard Gaussian variables. For every $(a_1,\ldots,a_k)\in [0,2)^k$ and $(b_1,\ldots, b_k)\in \Z^k$, the conditional law of 
\[ \left( \frac{h_{T'_n}(i)-(1-a_i/2) \ln n}{\sqrt{(1-a_i/4) \ln n}},\, i\in [k] \right),\]
given that $\deg_{T'_n} (i) \ge \floor{a_i\ln n} +b_i$ for all $i\in [k]$, converges to the law of $(N_i)_{i\in [k]}$.
\end{theorem}

The advantage of shifting towards this statement is that $T'_n$ may be analysed through a construction of Kingman's coalescent on $n$ elements \cite{AddarioEslava17,Eslava21}. 
\refT{thm:depths Kingman} is the first intermediate result for \refT{thm:joint conv}. 

The next proposition is the second intermediate result for \refT{thm:joint conv}.
For $i\in \Z$ and $B\in \cB_{\R}$ let 
\begin{align}
X_{i}^\n(B)&:= \#\{v\in [n]: \deg_{T_n}(v)=\floor{\log n}+i,\, (h(v)-\mu_n)/\sigma_n\in B  \}, \label{dfn:X} \\
X_{\ge i}^\n(B)&:= \#\{v\in [n]: \deg_{T_n}(v)\ge \floor{\log n}+i,\, (h(v)-\mu_n)/\sigma_n\in B\}. \label{dfn:X+}
\end{align}

\begin{proposition}\label{prop:moments}
Fix $c\in (0,2)$ and $K, M\in \N$. Let $(a_k)_{k\in [K]}$ be non-negative integers such that $
\sum_{k\in [K]} a_k=M$. Let $(j_k)_{k\in[K]}$ be a non-decreasing sequence of integers with $0\le K'=\min\{k:\, j_{k+1}=j_K\}$ and let $(B_k)_{k\in K}$ be sets in $\cB_\R$ satisfying $B_k\cap B_\ell=\emptyset$ whenever $j_k=j_\ell$ and $k\neq \ell$. If $0\le j_1+\log n< j_K+\log n<c\ln n$ then we have
\begin{align*}
&\E{\prod_{k=1}^{K'} \left(X^\n_{j_k}(B_k)\right)_{a_k}\prod_{k=K'+1}^{K} \left(X^\n_{\ge j_K}(B_k)\right)_{a_k}}\\
=& \prod_{k=1}^{K'} \left(2^{-j_k+\eps_n-1} \Phi(B_k)\right)^{a_k}
\prod_{k=K'+1}^{K} \left(2^{-j_K+\eps_n} \Phi(B_k)\right)^{a_k}
(1+o(1)).
\end{align*}
\end{proposition}

The proof of \refP{prop:moments} is given in \refS{sec:Proof moments}. We briefly explain an ingredient that is key in proving \refP{prop:moments} which is based on the exchangeability of vertices in the representation of $T_n$ as a Kingman's coalescent. As $T'_n$ is a relabelling of $T_n$, we have, jointly for all $i\in \Z$ and $j\in \N$,  
\begin{equation*}
\#\{v\in [n]: \deg_{T_n}(v)=i,\, h_{T_n}(v)=j \}\eqdist \#\{v\in [n]: \deg_{T'_n}(v)=i,\, h_{T'_n}(v)=j\}.
\end{equation*}

Note that each $X^\n_j(B)$ and $X^\n_{\ge j}(B)$ is a sum of indicator variables. 
Therefore, the factorial moments in \refP{prop:moments} are reduced to a sum of probabilities in terms of the degrees and depths of vertices in $T'_n$ which have the following form: \begin{align*}
&\p{\deg_{T'_n}(j)\ge m_j,\, h_{T'_n}(j)\in B_j,\, j\in [k]} \\
=&\p{\deg_{T'_n}(j)\ge m_j,\, j\in [k]}\p{h_{T'_n}(j)\in B_j,\, j\in [k]\,|\, \deg_{T'_n}(j)\ge m_j,\, j\in [k]};
\end{align*}
where $B_j\in \cB_\R$ and $m_j<2\ln n$ for each $j\in [k]$. The joint tails of finitely many vertices in $T'_n$ have been analyzed in \cite{AddarioEslava17} and we present the relevant result as \refP{prop:degrees Kingman}. The law of the depths conditional on $\{ \deg_{T'_n}(j)\ge m_j,\, j\in [k]\}$ is approximated using \refT{thm:depths Kingman}.

For \refT{thm:joint conv}, we next formulate the result as a statement about convergence of marked point processes. Let $\Z^*:=\Z\cup \{\infty\}$. Endow $\Z^*$ with the metric defined by $d(i,j):=|2^{-i}-2^{-j}|$ and $d(i,\infty):=2^{-i}$ for $i,j\in \Z$. Recall that $\cP$ is a Poisson point process on $\R$ with intensity rate $\lambda(x)=2^{-x}\ln 2$ and let $(\xi_x)_{x\in \cP}$ be an independent collection of i.i.d. standard Gaussian variables. 

Let $\eps\in [0,1]$. For $n\ge 1$, denote
\begin{align}\label{dfn:normalization}
\mu_n:=\left(1-\frac{\log e}{2}\right)\ln n\qquad \text{and} \qquad \sigma^2_n:=\left(1-\frac{\log e}{4}\right)\ln n.   
\end{align}

We define a ground process $\cP^\eps$ on $\Z^*$, and a marked point process $\cM\cP^\eps$ on $\Z^*\times \R$ given by  
\begin{align}
    \cP^\eps &:=\sum_{x\in \cP} \delta_{\floor{x+\eps}}, &     \cM\cP^\eps &:=\sum_{x\in \cP} \delta_{(\floor {x+\eps},\xi_x)}, \\ \intertext{Similarly, for all $n\in \N$ we use $T_n$ to define empirical point processes given by}
    \cP^\n &:=\sum_{v\in [n]} \delta_{\deg_{T_n}(v)-\floor{\log n}}, &    \cM\cP^\n &:=\sum_{v\in [n]} \delta_{(\deg_{T_n}(v)-\floor{\log n},(h(v)-\mu_n)/\sigma_n)}.
\end{align}

Let $\cM^{\#}_{\Z^*}$ and $\cM^{\#}_{\Z^*\times \R}$ be the space of boundedly finite measures of $\Z^*$ and $\Z^*\times \R$, respectively. It is clear that $\cM\cP^\eps$ and $\cM\cP^\n$ are elements of $\cM^{\#}_{\Z^*\times \R}$ (see, e.g. \cite[Definition~9.2.I]{DaleyVere-JonesII}).  The advantage of working on the state space $\Z^*$ is that intervals $[i, \infty]$ are compact and so the convergence of finite dimensional distributions of $\cM\cP^\nl$ implies, e.g., convergence in distribution of $X_{\ge i}^\nl(B)$ for any $B\in \cB_\R$ and $i\in \N$. We note that the weak convergence $\cP^\nl\to \cP^\eps$, for suitable sequences $(n_\ell)_{\ell\ge 1}$, had already been established in \cite{AddarioEslava17}.

Finally, the statement of \refT{thm:joint conv} is equivalent to the weak convergence of $\cM\cP^\nl\to \cM\cP^\eps$ in $\cM^{\#}_{\Z^*\times \R}$, for suitable sequences $(n_\ell)_{\ell\ge 1}$. The joint distribution of the variables in \refP{prop:moments} are finite dimensional distributions of the marked point process $\cM\cP^\nl$. Therefore, \refP{prop:moments} together with the method of moments establish the desired weak convergence. The complete details, together with the proof of \refC{cor:max conv}, are given in \refS{sec:Main proofs}.

\subsection{Kingman's coalescent approach for \refT{thm:depths Kingman}}\label{sec:Kingman}

For each $n\ge 1$, we consider the set of forests $\cF_n:=\{F:V(F)=[n]\}$ with vertex set $[n]$. 
An $n$-chain is a sequence $(F_n,\ldots,F_1)$ of elements of $\cF_n$ if for $1< i\le n$, $F_{i-1}$ is obtained from $F_i$ by adding an edge connecting two of the roots in $F_i$; note that $F_n$ is consequently the forest with $n$ one-vertex trees. We write $F_i=\{T_1^{(i)}, \ldots, T_{i}^{(i)}\}$, listing the trees according to an increasing order of their smallest-labeled vertex. 

\begin{figure} \begin{center}
 \includegraphics[scale=.75]{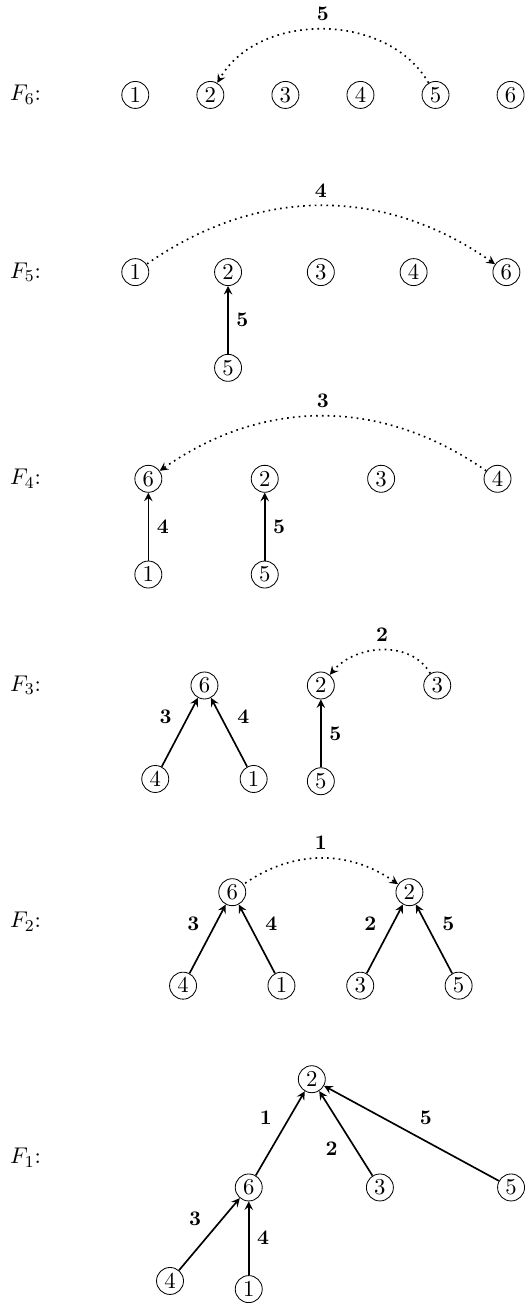}
 \end{center}
 \caption{An example of Kingman's $n$-coalescent $\mathbf{C}=(F_n,\ldots,F_1)$ for $n=6$. For $1< i\le n$, we represent the edge in $E(F_{i-1})\setminus E(F_i)$ with a dotted line in $F_i$. Edges are marked with the index of the first forest where $e$ is present.
 In this case, $\xi_6=\xi_4=\xi_3=1$, $\xi_5=\xi_2=0$ and $\{a_6,b_6\}=\{2,5\},$ $\;\{a_5,b_5\}=\{1,5\},\;\{a_4,b_4\}=\{1,4\},\; \{a_3,b_3\}=\{2,3\},\;\{a_2,b_2\}=\{1,2\}$.}
\label{fig:chain}
\end{figure}

We next define the construction of Kingman's $n$-coalescent as a random $n$-chain $\mathbf{C}=(F_n,\ldots, F_1)$. 
Informally, at each step $1\le j<n$, two uniformly random trees of $F_{n-j+1}$ are merged by adding an edge between their roots, the direction of the edge is determined independently and with equal probability. For an example of the process see \refF{fig:chain}.

\begin{definition}[Kingman's $n$-coalescent]\label{dfn:coalescent}
Let $n\in \N$. For each $1< i\le n$, choose $\{a_i,b_i\}\subset \{\{a,b\}: 1\le a<b\le i\}$ independently and uniformly at random; also let $(\xi_i)_{1< i\le n}$ be a sequence of independent $\mathrm{Bernoulli}(1/2)$ random variables.

The $n$-chain $\mathbf{C}=(F_n,\ldots, F_1)$ is defined as follows. For $1< i\le n$, starting from $F_i$, add an edge $e_{i-1}$ between the roots of $T_{a_{i}}^{(i)}$ and $T_{b_{i}}^{(i)}$; direct $e_{i-1}$ towards $r(T_{a_{i}}^{(i)})$ if $\xi_{i}=1$, and towards $r(T_{b_{i}}^{(i)})$ otherwise. Then $F_{i-1}$ contains the new tree and the remaining $i-2$ unaltered trees from $F_{i}$. 
\end{definition}

\begin{figure} \begin{center}
\includegraphics[scale=.9]{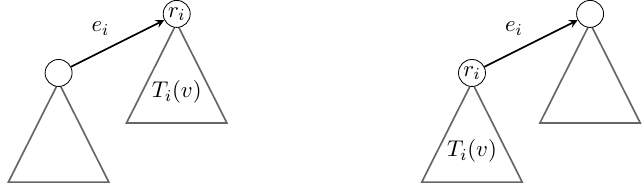}
\end{center}
\caption{For $v\in [n]$ and $1< i\le n$ let $r_i=r(T_i(v))$ and suppose $i\in \cS_n(v)$. Conditional on $F_i$, if $e_i$ is directed towards $r_i$, then the degree of $r_i$ increases by one in $F_{i-1}$. If $e_i$ is directed outwards $r_i$, then the depth of each $u\in T_i(v)$ increases by one in $F_{i-1}$.}
\label{fig:step t}
\end{figure}

It is not difficult to see that $\mathbf{C}$ has a uniform distribution on the $n!(n-1)!$ possible $n$-chains representing Kingman's $n$-coalescent, \cite{AddarioEslava17}. On the other hand, recall that $T_n$ has a uniform distribution on the $(n-1)!$ increasing trees in $\cI_n$. The next result essentially follows since bijections preserve the uniform measure on finite probability spaces; for details see \cite[Proposition 3.3]{AddarioEslava17} and also \cite[Propositions 1.2 and 1.3]{Eslava21}.

\begin{proposition}[Proposition 3.3 in \cite{AddarioEslava17}]
Let $n\ge 1$. Let $T'_n$ be defined as in \eqref{dfn:En} and let $T_1^{(1)}$ correspond to the final tree in $\mathbf{C}$ as constructed in \refD{dfn:coalescent}. Then $T'_n\eqdist T_1^{(1)}$; that is, for any labelled tree $T$ with vertex set $[n]$ we have
$\p{T'_n=T}=\p{T_1^{(1)}=T}$.
\end{proposition}

In the remainder of the paper, we analyse $T'_n$ using the construction in \refD{dfn:coalescent}; that is, we analyse the $n$-chain $(F_n,\ldots, F_1)$ and identify the final tree $T_1^{(1)}=T'_n$. To ease the notation, write 
\begin{align}
    d_n(v)=\deg_{T'_n}(v)\qquad \text{and} \qquad h_n(v)=h_{T'_n}(v),
\end{align}
for the degrees and depths of vertices in $T'_n$. 

For each vertex $v\in [n]$ and $1< i\le n$, let $T_i(v)$ be the tree in $F_i$ that contains $v$. Let $(s_{i,v})_{1<i\le n}$ and $(h_{i,v})_{1<i\le n}$ be sequences where $s_{i,v}$ is the indicator that $T_i(v)\in \{T_{a_i}^{(i)},T_{b_i}^{(i)}\}$ and $h_{i,v}$ is the indicator that the tail of the edge $e_{i-1}$ is $r(T_i(v))$. Note that $h_{i,v}=1$ only if $s_{i,v}=1$. Indeed, if $s_{i,v}=0$ then, when constructing $F_{i-1}$ from $F_i$, the trees that are merged do not include $T_i(v)$ and so necessarily $h_{i,v}=0$.

Now, conditional on $s_{i,v}=1$, there are two equally likely events regarding the merger of $T_i(v)$ to form $F_{i-1}$: if $e_{i-1}$ is directed towards $r(T_{i}(v))$, then the degree of $r(T_{i}(v))$ increases by one in $F_{i-1}$; otherwise, $e_{i-1}$ is directed out of $r(T_{i}(v))$ and all vertices in $T_{i}(v)$ increase their depth by one in $F_{i-1}$; see Figure~\ref{fig:step t}.

The choice of trees to be merged at each step is both independent and uniform so, for fixed $v\in [n]$, the variables $(s_{i,v})_{1< i\le n}$ are independent Bernoulli random variables with $\E{s_{i,v}}=2/i$. Moreover, since the variables $(\xi_i)_{i\in [n-1]}$ are independent from the selection of the trees to be merged, the variables $(h_{i,v})_{1< i\le n}$ are also independent Bernoulli random variables with $\E{h_{i,v}}=1/i$.

For $v\in [n]$, let 
\begin{align}
    \cS_n(v)&:=\{1< i\le n \,:\, s_{i,v}=1 \}, \label{dfn:S}\\
    \cS_{n,1}(v)&:=\{\ln^2 n< i\le n \,:\, s_{i,v}=1 \}; \label{dfn:S'}
\end{align}
we call these the \emph{selection} and \emph{truncated selection
sets} of $v$ and write $S_n(v):=|\cS_n(v)|$ and $S_{n,1}(v):=|\cS_{n,1}(v)|$, respectively. Using Lindeberg's condition (e.g., see \cite[Theorem 3.4.5]{Durret96}) we can see that, as $n \to \infty$, 
\begin{align}
\frac{S_n(v)-2\ln n}{\sqrt{2\ln n}}\todist N(0,1) 
\qquad \text{and} \qquad
\frac{S_{n,1}(v)-2\ln n}{\sqrt{2\ln n}}\todist N(0,1). \label{eq:S lims}
\end{align}

For each $1<i\le n$, the law of $h_{i,v}$ conditional on $s_{i,v}=1$ is Bernoulli  
with mean $1/2$. Since $(s_{i,v})_{1<i\le n}$ are independent, we have that 
\begin{align}\label{eq:h Ber sum}
    h_n(v)=\sum_{i=2}^n h_{i,v} \eqdist \sum_{j=1}^{S_n(v)} I_j,
\end{align}
where $(I_j)_{j\ge 1}$ are i.i.d.~Bernoulli random variables with mean $1/2$ independent of $\cS_n(v)$. Using the first expression of $h_n(v)$ in \eqref{eq:h Ber sum} and Lindeberg's condition, we have that $h_n(v)$ is asymptotically normal for $v\in [n]$ as $n\to \infty$. We may also establish the limiting distribution of $h_n(v)$ from the last expression in \eqref{eq:h Ber sum} and the following lemma; we include its short proof for completeness.  

\begin{lemma}\label{lem:skorohod}
Let $a\in [0,2)$, $b\in \Z$ and set $m:=\floor{a\ln n}+b$. Let $(Q_n)_{n\ge 1}$ be integer valued random variables such that the limiting law of $\frac{Q_n-2\ln n}{\sqrt{2\ln n}}$ is a standard Gaussian distribution. Let $N$ be a standard Gaussian variable and let $\hat{H}_n$ have the law\footnote{Abusing notation we let $\Bin{0,p}\equiv 0$.} of a $\Bin{Q_n-m,1/2}$ conditional on $|Q_n|\ge m$. Then,  
\begin{align*}
    \frac{\hat{H}_n-(1-a/2)\ln n}{\sqrt{(1-a/4)\ln n}}\todist N.
\end{align*}
\end{lemma}

\begin{proof}
Write $\hat{Q}_n$ for the law of $Q_n-m$ conditionally given that $Q_n\ge m$. Since $a<2$, we have that $\p{Q_n\ge m}\to 1$ as $n\to \infty$ and so the limiting law of $\frac{\hat{Q}_n-(2-a)\ln n}{\sqrt{2\ln n}}$ is a standard Gaussian distribution. Let $(I_1^{(n)},\ldots, I_n^{(n)})_{n\ge 1}$ be i.i.d.~Bernoulli random variables with mean $1/2$ independent of $(\hat{Q}_n)_{n\ge 1}$ and let $(X_n,Y_n)_{n\ge 1}$ be defined as
\begin{align*}
    X_n:=\frac{2\sum_{j=1}^n I_j^{(n)}-n}{\sqrt{n}} \qquad \text{and} \qquad Y_n:=\frac{\hat{Q}_n-(2-a)\ln n}{\sqrt{2\ln n}};
\end{align*}
by the independence assumption we have that $(X_n,Y_n)\todist (N_1,N_2)$ where $N_1$ and $N_2$ are independent standard Gaussian variables. 
By the Skorohod embedding theorem (e.g., see \cite[Theorem 6.7]{Billingsley}), we may couple $(\hat{Q}_n)_{n\ge 1}$ and $(I_1^{(n)},\ldots, I_n^{(n)})_{n\ge 1}$ in such a way that 
$(X_n,Y_n)\to (N_1,N_2)$ almost surely. 

Let $\hat{H}_n=\sum_{j=1}^{\hat{Q}_n} I_j^{(\hat{Q}_n)}$ and let $F(x)$ be the cumulative distribution function of  $\frac{\hat{H}_n-(1-a/2)\ln n}{\sqrt{(1-a/4)\ln n}}$. We have $\hat{Q}_n/\ln n\to 2-a$ and, in particular, $\hat{Q}_n\to \infty$, then for $x\in \R$ we have 
\begin{align*}
    F(x)=&\p{2\sum_{j=1}^{\hat{Q}_n}I_j^{(\hat{Q}_n)}\le x\sqrt{(4-a)\ln n}+(2-a)\ln n}\\
    =&\p{\frac{2\sum_{j=1}^{\hat{Q}_n}I_j^{(\hat{Q}_n)}
    -\hat{Q}_n}{\sqrt{\hat{Q}_n}} 
    \sqrt{\frac{\hat{Q}_n}{(4-a)\ln n}}
    +\frac{\hat{Q}_n - (2-a)\ln n}{\sqrt{2\ln n}}
    \sqrt{\frac{2}{(4-a)}}\le x}\\
    \to& \p{\sqrt{\frac{2-a}{4-a}}N_1 + \sqrt{\frac{2}{4-a}}N_2 \le x}=\Phi(x); 
\end{align*}
from which the desired convergence is established. 
\end{proof}

We will see in \refS{sec:main} that conditional on $d_n(v)\ge m$, we have $h_n(v) \eqdist \Bin{S_n(v)-m,1/2}$, which together with \refL{lem:skorohod}, yields the convergence in law of $h_n(v)$ conditional on $d_n(v)\ge m$ (assuming $2\ln n-m\to \infty$); for the details see \eqref{eq:hFj} and the discussion thereafter. Briefly sketched, this is the case $k=1$ of \refT{thm:depths Kingman}. 

For the general case, we have to deal with the correlations between the selection sets of finitely many vertices. For example, let $v$ and $w$ be distinct vertices in $T'_n$ and let $\lambda_{v,w}:=\max\{ 1< \ell\le  n:\, \ell \in \cS_n(v)\cap \cS_n(w)\}$; that is, in $F_{\lambda_{v,w}}$, vertices $v$  and $w$ belong to distinct trees that are merged in $F_{\lambda_{v,w}-1}$. The correlation between $(h_{i,v})_{2\le i\le n}$ and $(h_{i,w})_{2\le i\le n}$ is evident since  $h_{v,\lambda_{v,w}}=1-h_{w,\lambda_{v,w}}$ while $h_{i,v}=h_{i,w}$ for $i<\lambda_{v,w}$. In words, exactly one of $v$ or $w$ increases its depth in $F_{\lambda_{v,w}-1}$ and, from then on, their depths increase simultaneously for the remainder of the process.

To circumvent this problem, we analyse the coalescent at $F_{\floor{\ln^2 n}}$ where, for finitely many vertices, their selection sets --and so their depths-- are asymptotically independent (see \refL{lem:cond depth}) and the main contribution to their depths in $T'_n$ is already established (see \refL{lem:neglig depth}).

For each $v\in [n]$, define\footnote{If we were to define $h_{n,1}(i)=h_{F_j}(i)$ with $j=o(\sqrt{\ln n})$, then the statement of \refL{lem:neglig depth} would hold immediately; however, establishing \refL{lem:cond depth} would become a much more delicate task.}
\begin{align}\label{dfn:h1-2}
h_{n,1}(v):=h_{F_{\floor{\ln^2 n}}}(v)=\sum_{j\in \cS_{n,1}(v)} h_{j,v} 
\quad \text{and}\
\quad 
h_{n,2}(v):=h_n(v)-h_{n,1}(v).  
\end{align}

\begin{lemma}\label{lem:cond depth}
Let $k\in \N$. Let $(N_i)_{i\in [k]}$ be i.i.d.~standard Gaussian variables. For every $(a_1,\ldots, a_k)\in [0,2)^k$ and $(b_1,\ldots, b_k)\in \Z^k$, the conditional law of
\begin{align*}
    \left(\frac{h_{n,1}(i)-(1-a_i/2)\ln n}{\sqrt{(1-a_i/4)\ln n}},\, i\in [k]\right),
\end{align*}
given that $d_n(i)\ge \floor{a_i\ln n}+b_i$ for all $i\in [k]$, converges to the law of $(N_i)_{i\in [k]}$.
\end{lemma}

\begin{lemma}\label{lem:neglig depth}
Fix $k\in \N$ and $c\in (0,2)$. Let $m_i=m_i(n)<c\ln n$ for all $i\in[k]$. For each $\eps>0$ and any $j\in [k]$,
\[\p{h_{n,2}(j)\ge \eps\sqrt{\ln n} \,|\, d_n(i)
\ge m_i,\, i\in [k]} \to 0.\]
\end{lemma}

\refT{thm:depths Kingman} follows directly from Lemmas \refand{lem:cond depth}{lem:neglig depth} (e.g., see \cite[Theorem 3.1]{Billingsley}).

\section{Depths and degrees via truncated selection sets}\label{sec:main}

We briefly analyse the case of a single vertex. Write $\cS_n(1)=\{j_1,\ldots , j_{S_n(1)}\}$ with $j_1>j_2>\cdots >j_{S_n(1)}$. Throughout the coalescent process, both the degree and depth of 1 are determined by the variables $(h_{1,j_\ell})_{j_\ell\in \cS_n(1)}$ as follows.
Vertex 1 is the root of $T_n(1)$ in $F_n$. For $j\notin \cS_n(1)$, we have $T_{j}(1)=T_{j-1}(1)$. For $j\in \cS_n(1)$, if $h_{1,j}=0$ then the root of $T_{j}(1)$ increases its degree by one; otherwise, every vertex in $T_{j}(1)$ increases its depth by one.  In particular, the degree of 1 in $F_j$ corresponds to the first streak of values 0 in $(h_{v,j_\ell})_{j_\ell > j}$. More precisely, for $1\le j<n$,
$\deg_{F_j}(1)=\min\{d>j:\,h_{1,j_1}=h_{1,j_2}=\cdots =h_{1,j_d}=0\}$
and
\begin{align}\label{eq:hFj}
    h_{F_j}(1)=\sum_{\ell=1}^{S_n(1)\setminus [j]} h_{1,j_\ell}.
\end{align}

Recall that $d_n(1)=\deg_{F_1}(1)$ and $h_n(1)=h_{F_1}(1)$. Let $m\in \N$. The event $d_n(1)\ge m$ may occur only if $S_n(1)\ge m$ and, since $(h_{1,j_\ell})_{j_\ell\in \cS_n(1)}$ are i.i.d. Bernoulli variables with mean $1/2$, we have that $\p{d_n(1)\ge m|\, S_n(1)\ge m}=2^{-m}$. In addition, conditional on $d_n(1)\ge m$ we have $h_n(1) \eqdist \Bin{S_n(1)-m,1/2}$; this is well-defined since $d_n(1)\ge m$ implies $S_n(1)\ge m$. Together with \refL{lem:skorohod}, the distributional equivalence establishes the case $k=1$ of \refT{thm:depths Kingman}. 

The proof of \refL{lem:cond depth} adapts this argument to obtain the asymptotic limit of the variables $(h_{n,1}(i),\, i\in [k])$; that is, the depths of $k$ distinct vertices in $F_{\floor{\ln^2 n}}$. 
Recall the definition of $\cS_{n,1}(v)$ in \eqref{dfn:S'} and that $S_{n,1}(v)=|\cS_{n,1}(v)|$. Let 
\begin{align}
\Omega_1:=\{J\,:\, J\subset \{\floor{\ln^2 n}+1,\ldots, n\}\};  
\end{align}
we use a superscript bar to denote vectors indexed by $[k]$; for example, we write  $\bar{\cS}_{n,1}:=(\cS_{n,1}(1), \ldots, \cS_{n,1}(k))\in \Omega_1^k$. 

For $\bar{m}=(m_1,\ldots, m_k)\in \N^k$ and $\delta\in (0,2)$. Let 
\begin{align}
\cA_{\bar{m}}&:=
\{ \bar{J}\in \Omega^k_1:\, \p{\bar{\cS}_{n,1}=\bar{J},\, d_n(i)
\ge m_i,\, i\in [k]} >0\} \label{dfn:Am}
\end{align}
contain the outcomes of the first $k$ truncated sets $\bar{\cS}_{n,1}$ on the event $\{d_n(i) \ge m_i,\, i\in [k]\}$ and let
\begin{align}
\cB_\delta&:=
\{\bar{J}\in \Omega_1^k:\, (J_1,\ldots, J_k) \text{ are pairwise disjoint, } ||J_i|-2\ln n|\le\delta\ln n,\, i\in [k]\}. \label{dfn:Bnd}
\end{align}

The following lemma gathers two observations about events conditional on $\bar{\cS}_{n-1}=\bar{J}$ and a sufficient condition on $\bar{m}$ for $\cB_\delta\subset \cA_{\bar{m}}$ to hold.

\begin{lemma}\label{lem:cond J}
Let $k\ge 1$ and $\bar{m},\bar{\ell}\in \Z_{\ge 0}^k$. If $\bar{J}\in \Omega_1^k$ has pairwise disjoint sets and $|J_i|\ge m_i$ for $i\in [k]$, then
\begin{align}
&\p{h_{n,1}(i)\le \ell_i,\, d_n(i)\ge m_i;\, i\in [k]\,|\, \bar{\cS}_{n,1}=\bar{J}} \nonumber \\ =&
\prod_{i=1}^k \p{h_{n,1}(i)\le \ell_i,\, d_n(i)\ge m_i\,|\, \cS_{n,1}(i)=J_i};\label{eq:cond indep}
\end{align}
if $\bar{J}\in \cA_{\bar{m}}$ then
\begin{align}\label{eq:cond deg}
\p{d_n(i)\ge m_i;\, i\in [k]\,|\, \bar{\cS}_{n,1}=\bar{J}}=2^{-\sum_i m_i}.
\end{align}
Finally, if $\delta\in (0,2)$ satisfies $m_i<(2-\delta)\ln n$ for $i\in [k]$, then
    $\cB_\delta\subset \cA_{\bar{m}}$.
\end{lemma}

\begin{proof}
Let $k\ge 1$, $\bar{m}\in \Z_{\ge 0}^k$ and let $\bar{J}\in \Omega_1^k$ satisfy $|J_i|\ge m_i$ for $i\in [k]$. For each $i\in [k]$, write $J_i=\{j_{i,1},j_{i,2},\ldots,j_{i,|J_i|}\}$ with $j_{i,1}>j_{i,2}>\cdots>j_{i,|J_i|}$ and let $D_i=\{j_{i,1},j_{i,2},\ldots,j_{i,m_i}\}$. 

On the event that $\bar{\cS_{n,1}}=\bar{J}$, for each $i\in [k]$, 
\begin{align*}
h_{n,1}(i)=\sum_{\ell=1}^{|J_i|} h_{i,j_\ell},    
\end{align*}
and $d_n(i)\ge m_i$ if and only if $h_{i,j_\ell}=0$ for all $j_\ell\in D_i$. Consequently, for $\bar{\ell}\in \Z_{\ge 0}^k$, the event $\{h_{n,1}(i)\le \ell_i,\, d_n(i)\ge m_i,\,i\in [k]\}$ depends solely on the variables $(h_{i,j_{i,\ell}})_{j_{i,\ell}\in D_i}$ for $i\in [k]$. If $\bar{J}$ has pairwise disjoint sets, then so has $\bar{D}$ pairwise disjoint sets and \eqref{eq:cond indep} follows. 

Now, at each step of the coalescent process, exactly one vertex increases its degree. Observe that if $\bar{J}\in \cA_{\bar{m}}$ then $\p{d_n(i)\ge m_i,\,i\in [k]\,|\, \bar{\cS}_{n,1}=\bar{J}}>0$ and the indices in $\cup_{i\in[k]} D_i$ corresponds to steps at which precisely one of the vertices in $[k]$ increases its degree. It follows that $\bar{D}$ must have pairwise disjoint sets and so, conditional on $\bar{\cS}_{n,1}=\bar{J}$, the variables $(h_{i,j_{i,\ell}})_{j_\ell\in D_i}$, $i\in [k]$ are independent Bernoulli variables with mean $1/2$; establishing \eqref{eq:cond deg}.

Finally, if $\delta\in (0,2)$ satisfies $m_i<(2-\delta)\ln n$ for $i\in [k]$, then $\bar{J}\in \cB_\delta$ has pairwise disjoint sets and $|J_i|\ge m_i$ for each $i\in [k]$. Hence, $\p{d_n(i)\ge m_i\,|\, \cS_{n,1}(i)=J_i}>0$, and deterministically $h_n(i)\le n$ for all $i\in[k]$; we then infer from \eqref{eq:cond indep} with $\bar{\ell}=(n,\ldots, n)$ that $\bar{J}\in \cA_{\bar{m}}$. 
\end{proof}

The next lemma states that the distribution of $\bar{\cS}_{n,1}$ is asymptotically approximated by $k$ copies of independent truncated selection sets; we delay its proof to \refS{sec:selection sets}.

\begin{lemma}\label{lem:S bar}
Let $k\ge 2$. Let $\bar{\cR}_n:=(\cR_n(i),\, i\in [k])$ be $k$ independent copies of $\cS_{n,1}(1)$. For any $\delta\in (0,2)$, we have 
\begin{align}
\p{\bar{\cS}_{n,1}\in \cB_{\delta}}&=1+O(\ln^{-2} n), \label{eq:bulk S}
\intertext{and, uniformly for $\bar{J}\in \cB_{\delta}$, }
\p{\bar{\cS}_{n,1}=\bar{J}}&=(1+O(\ln^{-1} n))\p{\bar{\cR}_n=\bar{J}}. \label{eq:unif S}
\end{align}
\end{lemma}

Before proceeding to the proofs of Lemmas \ref{lem:cond depth}, \refand{lem:neglig depth}{lem:S bar} we recall a result from \cite{AddarioEslava17} which establishes asymptotic independence of the degrees of any finite number of vertices in $T'_n$.

\begin{proposition}[Proposition 4.2 in \cite{AddarioEslava17}] \label{prop:degrees Kingman}
Fix $c\in (0,2)$ and $k\in \N$. There exists $\beta=\beta(c,k)>0$ such that uniformly over positive integers $m_1,\ldots, m_k< c\ln n$,
\[\p{d_n(i)\ge m_i,\, i\in [k]}= 2^{-\sum_{i} m_i} (1+o(n^{-\beta})).\]
\end{proposition}

\subsection{Proof of \refL{lem:cond depth} }

Fix $k\in \N$, $\bar{a}\in [0,2)^k$ and $\bar{b}\in \Z^k$. Let $\bar{x}\in \R^k$, $m_i:=\floor{a_i \ln n}+b_i$ and $\ell_i:=(1-a_i/2)\ln n + x_i\sqrt{(1-a_i/4)\ln n}$; let $c\in (0,2)$ satisfy $m_i<c\ln n$ for all $i\in [k]$ and $n$ be large enough.

Let us write $m=m_1$. We infer from \eqref{eq:hFj} that, conditional on $\{d_n(1)\ge m, S_{n,1}(1)\ge m$\}, we have $h_{n,1}(1) \eqdist \Bin{\hat{S}_n,1/2}$; where $\hat{S}_n$ has the law of $S_{n,1}(1)-m$ conditional on $S_{n,1}(1)\ge m$. In particular,
\begin{align*}
    &\p{h_{n,1}(1)\le \ell_1 \,|\, d_n(1)\ge m\,, S_{n,1}(1)\ge m}\\
    =&\p{\Bin{S_{n,1}(1)-m,1/2}\le \ell_1\,|\, S_{n,1}(1)\ge m};
\end{align*}
from which we obtain  
\begin{align*}
    &\p{h_{n,1}(1)\le \ell_1 ,\, d_n(1)\ge m}\\
    =&\p{\Bin{S_{n,1}(1)-m,1/2}\le \ell_1\,|\, S_{n,1}(1)\ge m}\p{d_n(1)\ge m,\, S_{n,1}(1)\ge m} \\
    &+ \p{h_{n,1}(1)\le \ell_1 ,\, d_n(1)\ge m,\, S_{n,1}(1)< m}.
\end{align*}

Now, recall that $d_n(1)\ge m$ only if $S_n(1)\ge m$ and observe that the event $\{d_n(1)\ge m\}$ implies that  throughout the complete coalescent process there are at least $m$ fair, independent coin flips that land heads which are also independent from the selection sets. Since both $S_n(1)$ and $S_{n,1}(1)$ concentrate around $2\ln n$ we have 
\begin{align*}
    \lim_{n\to \infty} \frac{\p{d_n(1)\ge m,\, S_{n,1}(1)\ge m}}{\p{d_n(1)\ge m}}
    =\lim_{n\to \infty} \frac{2^{-m}\p{S_{n,1}(1)\ge m}}{2^{-m}\p{S_n(1)\ge m}} =1,
\end{align*}
and 
\begin{align*}
    \lim_{n\to \infty} \frac{\p{d_n(1)\ge m,\, S_{n,1}(1)< m}}{\p{d_n(1)\ge m}}
    \le \lim_{n\to \infty} \frac{2^{-m}\p{S_{n,1}(1)< m}}{2^{-m}\p{S_n(1)\ge m}} =0.
\end{align*}
Putting together all the estimates above, we obtain as $n\to \infty$
\begin{align}
   \p{h_{n,1}(1)\le \ell_1 \,|\, d_n(1)\ge m}\to  \p{\Bin{S_{n,1}(1)-m,1/2}\le \ell_1\,|\, S_{n,1}(1)\ge m};
\end{align}
which establishes, by \refL{lem:skorohod}, the desired convergence if $k=1$ and more generally, by the exchangeability of the vertices, for any $i\in [k]$,
\begin{align}\label{eq:lim i}
    \lim_{n\to \infty} \p{h_{n,1}(i)\le \ell_i\,|\, d_n(i)\ge m_i}= \Phi(x_i).
\end{align}

For $k\ge 2$, assume that  
\begin{align}
    &\p{h_{n,1}(i)\le \ell_i,d_n(i)\ge m_i,\,i\in [k]}\nonumber \\=&\prod_{i=1}^k\p{h_{n,1}(i)\le \ell_i,d_n(i)\ge m_i}+o\left(2^{-\sum_i m_i}\right). \label{eq:prod k}
\end{align}
By \refP{prop:degrees Kingman}, we have
\begin{align}
 \prod_{i=1}^k \p{d_n(i)\ge m_i}=(1+o(1))2^{-\sum_i m_i}=\p{d_n(i)\ge m_i,\,i\in [k]}.\label{eq:prod degrees}
\end{align} 
Putting together \eqref{eq:lim i}--\eqref{eq:prod degrees}  implies, for arbitrary $\bar{x}\in \R^k$, 
\begin{align*}
&\p{h_{n,1}(i)\le \ell_i,d_n(i)\ge m_i,\,i\in [k]}   \\ =&(1+o(1))\p{d_n(i)\ge m_i,\,i\in [k]}\prod_{i=1}^k\p{h_{n,1}(i)\le \ell_i\,|\,d_n(i)\ge m_i}\\
    =&(1+o(1))\p{d_n(i)\ge m_i,\,i\in [k]}\prod_{i=1}^k\Phi(x_i);
\end{align*}
this yields the desired convergence in law for $k\ge 2$. 

It then remains to prove \eqref{eq:prod k}. \refL{lem:cond J} and the choice of $c$ gives $\cB_{2-c}\subset \cA_{\bar{m}}$, together with \refL{lem:S bar} we have
\begin{align}\label{eq:negl}
    \p{\bar{\cS}_{n,1}\in \cA_{\bar{m}}\setminus \cB_{n,2-c}}\le \p{\bar{\cS}_{n,1}\in \Omega_1^k\setminus \cB_{n,2-c}}=o(1). 
\end{align}
Now, we decompose the event in the left-hand side of \eqref{eq:prod k} and use \eqref{eq:cond indep}, \eqref{eq:cond deg} and \eqref{eq:negl} to neglect the contribution of $\bar{J}\notin \cB_{2-c}$,
\begin{align*}
&\p{h_{n,1}(i)\le \ell_i,d_n(i)\ge m_i,\,i\in [k]}
\\
=&\sum_{\bar{J}\in \cB_{n,2-c}} \hspace{-.25cm} \p{h_{n,1}(i)\le \ell_i,d_n(i)\ge m_i,\, i\in [k]\,|\, \bar{\cS}_{n,1}=\bar{J}}\p{\bar{\cS}_{n,1}=\bar{J}} +o\left(2^{-\sum_i m_i}\right)\\
=&\sum_{\bar{J}\in \cB_{n,2-c}} \prod_{i=1}^k
    \p{h_{n,1}(i)\le \ell_i,d_n(i)\ge m_i\,|\, \bar{\cS}_{n,1}=\bar{J}}\p{\bar{\cS}_{n,1}=\bar{J}} +o\left(2^{-\sum_i m_i}\right)\\
=&\sum_{\bar{J}\in \Omega_1^k} \prod_{i=1}^k
    \p{h_{n,1}(i)\le \ell_i,d_n(i)\ge m_i\,|\, \bar{\cS}_{n,1}=\bar{J}}\p{\bar{\cS}_{n,1}=\bar{J}} +o\left(2^{-\sum_i m_i}\right);
\end{align*}
establishing \eqref{eq:prod k}.

\subsection{Proof of \refL{lem:neglig depth}}

Let $k\ge 2$, $c\in (0,2)$ and $\eps>0$. Let $m_i=m_i(n)<c\ln n$ for all $i\in [k]$. By \refP{prop:degrees Kingman} and the exchangeability of the vertices, it suffices to show that
\begin{align*}
    \p{h_{n,2}(1)\ge \eps\sqrt{\ln n},\,d_n(i)\ge m_i,\, i\in[k]}=o\left(2^{-\sum_i m_i}\right).
\end{align*}
Let us assume that uniformly over $\bar{J}\in \cB_{n,2-c}$, 
\begin{align}\label{eq:neglig}
\p{h_{n,2}(1)\ge \eps\sqrt{\ln n},\,  d_n(i)
\ge m_i,\, i\in [k]\,|\, \bar{\cS}_{n,1}=\bar{J}}
= o\left(2^{-\sum_i m_i}\right). 
\end{align}
Then using \eqref{eq:neglig}, \eqref{eq:cond deg} and \eqref{eq:bulk S}, we get
\begin{align*}
&\quad \p{h_{n,2}(1)\ge \eps\sqrt{\ln n},\, d_n(i)
\ge m_i,\, i\in [k]}\\
=&\sum_{\bar{J}\in \cA_{\bar{m}}} \p{h_{n,2}(1)\ge \eps\sqrt{\ln n},\, d_n(i)\ge m_i,\, i\in [k] \;|\;
\bar{\cS}_{n,1}=\bar{J}} \p{\bar{\cS}_{n,1}=\bar{J}}\\
\le& \p{\bar{\cS}_{n,1}\in \cB_{2-c}}o\left(2^{-\sum_i m_i}\right) + 2^{-\sum_i m_i}\p{\bar{\cS}_{n,1}\in\cA_{\bar{m}}\setminus \cB_{2-c}})= o\left(2^{-\sum_i m_i}\right).
\end{align*}

Now, to prove \eqref{eq:neglig}, let $\bar{J}\in \cB_{2-c}$ and simplifying notation write $[\ln^2 n]=[\floor{\ln^2 n}]$. Observe that $|J_i|\ge m_i$ for $i\in [k]$ and so the selection times $\cS_n(1)\cap [\ln^2 n]$ are independent of the event $\{d_n(i)\ge m_i,\,i\in [k]\}\cap \{\bar{\cS}_{n,1}=\bar{J}\}$. Moreover, $|\cS_n(1)\cap [\ln^2 n]|\eqdist S_{\floor{\ln^2 n}+1}$ which has mean $(1+o(1))2\ln\ln n$. Using Berstein's inequalities (e.g., see \cite[Theorem 2.10]{BoucheronLugosiMassart}) we get  
\begin{align*}
&\p{|\cS_n(1)\cap [\ln^2 n]|\ge \eps\sqrt{\ln n} \; \middle| \; d_n(i)
\ge m_i,\, i\in [k],\, \bar{\cS}_{n,1}=\bar{J}}\\
=&\p{S_{\floor{\ln^2 n}+1}(1)\ge \eps\sqrt{\ln n}}=o(1).
\end{align*}
It is straightforward that $h_{n,2}(1)$ is stochastically dominated by $|\cS_n(1)\cap [\ln^2 n]|$, so the estimate above together with \eqref{eq:cond deg} (recall that $\cB_{2-c}\subset \cA_{\bar{m}}$) establishes \eqref{eq:neglig}.

\subsection{Proof of \refL{lem:S bar}}\label{sec:selection sets}

Let $k\ge 2$ and $$\tau_k:=\max\{1< j\le n: s_{j,i}=s_{j,i'}=1 \text{ for distinct } i,i'\in [k]\}.$$ Then,  
\begin{align}
\p{(\cS_{n,1}(i),\, i\in [k]) \text{ are pairwise disjoint}}&=\p{\tau_k\le \ln^2 n }\nonumber\\
&\ge 1-2k^2\ln^{-2} n;\label{eq:tau}
\end{align}
the equality holds by definition of $\tau_k$ and the condition on  $\bar{\cS}_{n,1}$ having pairwise disjoint sets, while the inequality follows from \cite[Lemma 4.6]{AddarioEslava17}.

Let $\delta\in (0,2)$, by Berstein's inequality (e.g., see \cite[Theorem 2.10]{BoucheronLugosiMassart}) there is $\beta=\beta(\delta)>0$ such that for any $i\in [k]$
\begin{align}\label{eq:Berstein}
\p{|S_{n,1}(i)-2\ln n|\ge \delta \ln n}= o(n^{-\beta});    
\end{align}
it follows that a union bound, together with \eqref{eq:tau} and \eqref{eq:Berstein} yields \eqref{eq:bulk S}. 

For the proof of \eqref{eq:unif S} we will use the following factors. For $\ln^2n< m\le n$ let 
\begin{align}
p_{m,0}&:=\frac{(m-k)(m-k-1)}{m(m-1)}, && q_{m,0}:=\left(1-\frac{2}{m}\right)^k, \label{dfn:pq0}\\
p_{m,1}&:=\frac{2(m-k)}{m(m-1)},  && q_{m,1}:=\frac{2}{m}\left(1-\frac{2}{m}\right)^{k-1}. \label{dfn:pq1}
\end{align}
There exists a constant $C=C(k)>0$ such that for $m$ large enough, we have for $\sigma\in \{0,1\}$,
\begin{align}\label{eq:bound pq}
\left(1-\frac{C}{m^2}\I{\sigma=0}\right)q_{m,\sigma} & \;< \;p_{m,\sigma} \;< \;q_{m,\sigma} \left(1+\frac{C}{m}\I{\sigma=1}\right).
\end{align}
Indeed, write $q_{m,0}=1-\frac{2k}{m}+\frac{2k(k-1)}{m^2}+cm^{-3}$ and $q_{m,1}=\frac{2}{m}\left(1-\frac{2(k-1)}{m}+c'm^{-2}\right)$, it is straightforward that $\max\{|c|,|c'|\}\le 2^k$. We then have 
\begin{align*}
q_{m,0}-p_{m,0}&=\frac{k(k-1)(m-2)}{m^2(m-1)} +cm^{-3},\\
p_{m,1}-q_{m,1}&=\frac{2}{m} \left(\frac{(k-1)(m-2)}{m(m-1)} -c'm^{-2}\right);
\end{align*}
which implies that $q_{m,0}-p_{m,0}=\Theta(m^{-2})$ and $q_{m,1}-p_{m,1}=\Theta(m^{-2})$. 
Similarly, since $\lim_{m\to \infty} q_{m,0}=1$ and $\lim_{m\to \infty} mq_{m,1}=2$, there is $C=C(k)>0$ such that for $m$ large enough,
$0<q_{m,0}-p_{m,0}\le Cq_{m,0}m^{-2}$ and $0<p_{m,1}-q_{m,1}\le Cq_{m,1}m^{-1}$. 

In what follows fix $\bar{J}\in \cB_{\delta}$ and 
let $\bar{\cR}_n$ be $k$ i.i.d.~copies of $\cS_{n,1}(1)$. Set $r_{m,i}=\I{m\in \cR_n(i)}$ and $j_{m,i}=\I{m\in J_i}$ for $\ln^2n < m\le n$ and $i\in [k]$. Let $\sigma_m=\sum_{i\in [k]} j_{m,i}$ and observe that $\sigma_m\in \{0,1\}$. Recall that $\cS_{n,1}(i)=\{\ln^2 n< m\le n:\, s_{m,i}=1\}$ for $i\in[k]$ and let $A_m:=\{s_{m,i}=j_{m,i},\, i\in [k]\}$ for $\ln^2n< m\le n$. Then
\begin{align*}
\p{\bar{\cS}_{n,1}=\bar{J}}&=\p{A_n}\prod_{m=\floor{\ln^2 n}+1}^{n-1} \p{A_m\,|\, A_l,\,m<l\le n}=\prod_{m=\floor{\ln^2}+1 n}^n p_{m,\sigma_m};
\end{align*}
to see this,use that $\bar{J}$ has pairwise disjoint sets to infer that $\#\{T_m(i);\, i\in [k]\}=k$ for all $m> \ln^2 n$; that is, at no point we have two trees $T_m(i)$ and $T_m(i')$ selected to be merged with $i,i'\in [k]$ and $m> \ln^2 n$. Then, for $\sigma\in \{0,1\}$, $p_{m,\sigma}$ is the probability of selecting $\sigma$ trees from $\{T_m(i);\, i\in [k]\}$ when selecting two distinct trees uniformly at random from the forest $F_m$. 

Next, observe that, for $n$ large enough,
\begin{align*}
    \prod_{m=\floor{\ln^2 n}+1}^n \left(1-\frac{C}{m^2} \right) 
>1- \hspace{-.15cm}\sum_{m=\floor{\ln^2 n}+1}^n \frac{2C}{m^2} >1 -2C \int_{\ln^2 n}^{\infty} x^{-2}dx= 1-2C\ln^{-2} n.
\end{align*}
On the other hand, we have that $\sum \sigma_m= \sum_{i\in [k]} |J_i| \le  (2+\delta)k\ln n$ and, for $n$ large enough, 
\begin{align*}
    \prod_{m:\sigma_m=1} \left(1+\frac{C}{m}\right)\le \exp\left(\frac{(2+\delta)Ck\ln n}{\ln^2 n}\right)
<1+\frac{(2+\delta)2Ck}{\ln n};
\end{align*}
putting together these estimates we obtain
\begin{align*}
\prod_{m=\floor{\ln^2 n}+1}^n p_{m,\sigma_m}=(1+O(\ln^{-1} n))\prod_{m=\floor{\ln^2 n}+1}^n q_{m,\sigma_m}.
\end{align*}
The independence of the sets in $\bar{\cR}_n$ gives $\p{\bar{\cR}_n=\bar{J}}=\prod_{m=\floor{\ln^2 n}+1}^n q_{m,\sigma_m}$ and so $\p{\bar{\cS}_{n,1}=\bar{J}}=(1+O(\ln^{-1} n))\p{\bar{\cR}_n=\bar{J}}$; establishing \eqref{eq:unif S}.

\section{ Proof of \refP{prop:moments}}\label{sec:Proof moments}

Recall that we write $d_n(v)=\deg_{T'_n}(v)$ and $h_n(v)=h_{T'_n}(v)$ and, from \eqref{dfn:normalization}, 
$\mu_n=(1-(\log e)/2)\ln n$ and $\sigma_n=(1-(\log e)/4)\ln n$. A simpler version of the next lemma appeared in \cite[Lemma~5.1]{AddarioEslava17}. 

\begin{lemma}\label{lem:moments}
Let $k\in \N$. Let $(c_i)_{i\in [k]}$ be a non-decreasing sequence of integers and let $(B_i)_{i\in [k]}$ be sets in $\cB_\R$ such that $B_i\cap B_j=\emptyset$ whenever $c_i=c_j$ and $i\neq j$. Set $k'=\min\{i\ge 0:\, c_{i+1}=c_k\}$ and $m_i:=\floor{\log n}+c_i$ for $i\in [k]$, and
\begin{align*}
\cD_{\bar{m}} &:=
\{d_n(i)=m_i,d_n(j)\ge m_j\, 1\le i\le k'<j\le k\},\\
\cH_{\bar{B}} &:=
\left\{\frac{h_n(i)-\mu_n \ln n}{\sqrt{\sigma_n^2\ln n}}\in B_i,\, i\in [k] \right\}.
\end{align*}
Then
\[\p{\cD_{\bar{m}},\, \cH_{\bar{B}}}= (1+o(1)) 2^{-k'-\sum_i m_i} \prod_{i=1}^{k} \Phi(B_i).\]
\end{lemma}

\begin{proof}
By the inclusion-exclusion principle, 
\begin{align}\label{eq:inclusion}
\p{d_n(i)=m_i,\, i\in [k']}
= \sum_{\ell=0}^{k'} \sum_{\stackrel{S\subset [k']}{|S|=\ell}} (-1)^\ell
\p{d_n(i)\ge m_i+\I{i\in S},\, i\in [k']}.
\end{align}
By intersecting the event $\cH_{\bar{B}}\cap \{d_n(j)\ge m_j,\, k'<j\le k\}$ with $\{d_n(i)= m_i,\, i\in [k']\}$ and each of the events $\{d_n(i)\ge m_i+\I{i\in S},\, i\in [k']\}$ from \eqref{eq:inclusion} we obtain
\begin{align}\label{eq:inclusion2}
\p{\cD_{\bar{m}},\, \cH_{\bar{B}}}
=\sum_{\ell=0}^{k'} \sum_{\stackrel{S\subset [k']}{|S|=\ell}} (-1)^\ell 
\p{\cH_{\bar{B}}, d_n(i)\ge m_i+\I{i\in S},\, i\in [k]}.
\end{align}
By \refP{prop:degrees Kingman}, for each $S\subset [k']$ with $|S|=\ell$ we have $$\p{d_n(i)\ge m_i+\I{i\in S},\, i\in [k]}=2^{-\ell-\sum_i m_i}(1+o(1)).$$ 
Conditioning on $\{d_n(i)\ge m_i+\I{i\in S},\, i\in [k]\}$, by \refT{thm:depths Kingman} (with $a_i=\log e$ and $b_i=c_i+\I{i\in S}$), we get
\begin{align*}
\p{\cH_{\bar{B}},\, d_n(i)\ge m_i+\I{i\in S},\, i\in [k]}
=& (1+o(1)) 2^{-\ell-\sum_i m_i} \prod_{i=1}^{k} \Phi(B_i).
\end{align*}
Finally, we use the above estimate for each of the summands of the right-hand side of \eqref{eq:inclusion2} to obtain the desired expression; namely,
\begin{align*}
\p{\cD_{\bar{m}},\, \cH_{\bar{B}}}
=&(1+o(1))2^{-\sum_i m_i}\prod_{i=1}^{k} \Phi(B_i)\sum_{\ell=0}^{k'}\sum_{\stackrel{S\subset [k']}{|S|=\ell}} (-1)^\ell 2^{-\ell}\\
=&(1+o(1))2^{-k'-\sum_i m_i}\prod_{i=1}^{k} \Phi(B_i).
\end{align*}

\vspace{-.5cm}
\end{proof}

\begin{proof}[Proof of \refP{prop:moments}]
Fix $c\in (0,2)$ and $K, M\in \N$. Let $(a_k)_{k\in [K]}$ be non-negative integers such that $\sum_{k\in [K]} a_k=M$. Let $(j_k)_{k\in [K]}$ be a non-decreasing sequence of integers with $K'=\min\{k\ge 0:\, j_{k+1}=j_K\}$ such that $0\le j_1+\log n$ and $j_K+\log n <c\ln n$. Let $(B_k)_{k\in [K]}$ be sets in $\cB_\R$ satisfying $B_k\cap B_\ell=\emptyset$ whenever $j_k=j_\ell$ and $k\neq \ell$.  

We define $m_i\in \N$ and $A_i\subset \R$ as follows. For each $k\in [K]$, if $\sum_{\ell=1}^{k-1} a_\ell< i \le \sum_{\ell=1}^k a_\ell$ then set $m_i=\floor{\log n}+j_k$ and let $A_i=B_{j_k}$. Let $M'=\sum_{k\in [K']}a_k$. With foresight, we verify that $\prod_{i=1}^{M} \Phi(A_i)=\prod_{k=1}^{K} \Phi(B_k)^{a_k}$ and similarly,
\begin{align*}
M\log n-M'-\sum_{i=1}^M m_i &=\sum_{k=1}^{K'} (-j_k-1+\eps_n)a_k+ \sum_{k=K'+1}^{K}(-j_K+\eps_n)a_{k}.
\end{align*}

Now, consider the sets 
\begin{align*}
\cD_{\bar{m}} &:=\{d_n(i)=m_i,\, d_n(j)\ge m_j,\,1\le i\le M'<j\le M\},\\
\cH_{\bar{A}} &:=\left\{\frac{h_n(i)-\mu_n \ln n}{\sqrt{\sigma_n^2\ln n}}\in A_i,\, i\in [M] \right\}.
\end{align*}

Recall that $X_i^\n(B), X_{\ge i}^\n(B)$, $i\in \Z,B\in \cB_\R$ are sums of indicator functions; by the definition of $T'_n$ (as a relabelling of $T_n$) and the exchangeability of its vertex set, we have
\begin{align*}
\E{\prod_{k=1}^{K'} \left(X^\n_{j_k}(B_k)\right)_{a_k}\prod_{k=K'+1}^{K} \left(X^\n_{\ge j}(B_k)\right)_{a_k}}
=(n)_M \p{\cD_{\bar{m}},\,\cH_{\bar{A}}}.
\end{align*}
Finally, using the above estimates and since $(n)_M=n^M(1+o(n^{-1}))$, 
\refL{lem:moments} implies 

\begin{align*}
&\E{\prod_{k=1}^{K'} \left(X^\n_{j_k}(B_k)\right)_{a_k}\prod_{k=K'+1}^{K} \left(X^\n_{\ge j_k}(B_k)\right)_{a_k}}\\
=&(1+o(1))\left( 2^{M\log n-M'-\sum_i m_i} \right) \prod_{i=1}^{M} \Phi(A_i)\\
=&(1+o(1))\prod_{k=1}^{K'} \left( 2^{-j_k-1+\eps_n} \Phi(B_k)\right )^{a_k} \prod_{k=K'+1}^{K} \left( 2^{-j'-\eps_n} \Phi(B_k)\right)^{a_k},
\end{align*}
as desired.
\end{proof}

\section{Proof of \refT{thm:joint conv} and  \refC{cor:max conv}}\label{sec:Main proofs} 

Throughout the section we let $\eps\in [0,1)$ and consider an increasing subsequence $(n_\ell)_{\ell \ge 1}$ with $\log n_\ell-\floor{\log n_\ell}\to \eps$ as $\ell\to \infty$. Recall also that $\mu_{n_\ell}= (1-(\log e)/2) \ln n_\ell $ and $\sigma^2_{n_\ell}= (1-(\log e)/4) \ln n_\ell$ as in \eqref{dfn:normalization}. To simplify notation, let $\lambda_j:=2^{-j-1+\eps}$ for each $j\in \Z$. 

\begin{proof}[Proof of \refT{thm:joint conv}]
We will prove that  $\cM\cP^\nl$ converges weakly to $\cM\cP^\eps$, as $\ell \to \infty$, in $\cM^{\#}_{\Z^*\times \R}$. The weak convergence is equivalent to the statement in \refT{thm:joint conv} since we may list the vertices in $T_n$ in decreasing order of degree as $v^{(i)}$ and let $\cP=\{P_i,\, i\ge 1\}$ as defined in \eqref{dfn:P}; then
\begin{align*}
    \cM\cP^\eps =\sum_{i\ge 1} \delta_{(\floor{P_i}+\eps,\xi_{P_i})}, \quad \text{and}\quad 
    \cM\cP^\n =\sum_{i=1}^n \delta_{(\deg_{T_n}(v^{(i)})-\floor{\log n},(h(v^{(i)})-\mu_n)/\sigma_n)}.
\end{align*}

Consider the semiring of bounded sets that generates $\Z^*\times \R$ defined by
\begin{align}
\cA:=\{\{i\}\times (a,b],\, i\in \Z,\, a,b\in \R\}\cup \{[i,\infty]\times (a,b],\, i\in \Z,\, a,b\in \R\}.
\end{align}
It suffices to prove convergence of the counting measures of finite collections of disjoint sets in $\cA$; see \cite[Proposition~9.2.III and Theorem~11.1.VII]{DaleyVere-JonesII}. 

Consider any non-decreasing sequence of integers $(j_k)_{k\in[K]}$ with $0\le K'=\min\{k:\, j_{k+1}=j_K\}$ and let $(B_k)_{k\in K}$ with $B_k=(a_k,b_k]$ satisfying $B_k\cap B_\ell=\emptyset$ whenever $j_k=j_\ell$ and $k\neq \ell$.
Observe that the variables 
\begin{align}
    &(\cP^\eps(\{j_1\}\times B_1), \ldots, \cP^\eps(\{j_{K'}\}\times B_{K'}) \quad \text{and}\nonumber \\ &\cP^\eps([j_{K'+1},\infty]\times B_{K'+1}), \ldots,  \cP^\eps([j_{K},\infty]\times B_{K})) \label{eq:joint},
\end{align} 
are independent with distributions 
\begin{align*}
    \cP^\eps(\{i\},B)\eqdist \Poi{2^{-i-1+\eps}\Phi(B)} \qquad\text{and}\qquad
    \cP^\eps([i,\infty],B)\eqdist \Poi{2^{-i+\eps}\Phi(B)}.
\end{align*}
Recall that, if $X\eqdist \Poi{\lambda}$, then $\E{(X)_a}=\lambda^a$ for all integers $a\ge 0$. Setting $c:=3/2$, we infer that for $n$ large enough $0\le j_1+\log n< j_K+\log n<c\ln n$. It follows from \refP{prop:moments} and the method of moments  (e.g., see \cite[Section 6.1]{JLR})
 that 
$(X^\nl_{j_1}(B_1), \ldots, X^\nl_{j_{K'}}(B_{K'}), X^\nl_{\ge j_{K'+1}}(B_{K'+1}), \ldots,  X^\nl_{\ge j_{K}}(B_{K}))$ converge in distribution to the variables in \eqref{eq:joint};
completing the proof of the weak convergence of $\cM\cP^\nl$ to $\cM\cP^\eps$ in $\cM^{\#}_{\Z^*\times \R}$. 
\end{proof}

\begin{lemma}\label{lem:dist M}
Let $\eps>0$. The random variable $\mathbf M_\eps$ from \eqref{dfn:Me} is well defined. 
\end{lemma}

\begin{proof}
Let $\eps>0$. For each $m\in \Z$, let $p_m=(1-e^{-\lambda_m})e^{-\lambda_m}$. Observe that $2\lambda_m=\lambda_{m-1}$ and so 
\begin{equation*}
    \sum_{m\in \Z} p_m=\sum_{m\in \Z} e^{-\lambda_m}-e^{-\lambda_{m-1}} =\lim_{m\to \infty} e^{-\lambda_m}-\lim_{m\to -\infty} e^{-\lambda_{m-1}}=1.
\end{equation*}
Let $(Y_m,\, m\in \Z)$ be a collection of independent random variables where each $Y_m$ has the distribution of a Poisson random variable with mean $\lambda_m$ conditioned to being positive. Let $L_\eps$ be a random variable, independent of $(Y_m,\, m\in \Z)$ such that $\p{L_\eps=m}=p_m$, for $m\in \Z$. 

We will show that letting $\mathbf M_\eps=Y_{L_\eps}$ recovers the probability mass function in \eqref{dfn:Me}. To see this, observe that 
\begin{align*}
\p{\mathbf M_{\eps}=k} =\sum_{m\in \Z} \p{L_\eps=m}\p{Y_m=k}   =\sum_{m\in Z} (1-e^{-\lambda_m})e^{-\lambda_m} \cdot \frac{e^{-\lambda_m} \lambda_m^k}{k!(1-e^{-\lambda_m})},
\end{align*}
from which the lemma follows. 
\end{proof}

\begin{proof}[Proof of \refC{cor:max conv}]

Let $k\ge 1$. Since $\cP$ is independent of $(N_i)_{i\ge 1}$ it follows from \refT{thm:joint conv} that 
\begin{align*}
    \left((h_{T_{n_\ell}}(v^{(i)}) -\mu_{n_\ell})/\sigma_{n_\ell}
    ,\, i\in [k] \right) \todist \left(N_i,\, i\in [k]\right);
\end{align*}
on the other hand, conditional on $|\cM_{n_\ell}|=k$ we have $\cM_{n_\ell}=\{v^{(1)},\ldots,v^{(k)}\}$ and, \refT{thm:joint conv} gives 
\begin{align*}
    \left((h_{T_{n_\ell}}(v_i) - \mu_{n_\ell})/\sigma_{n\ell},\, v_i\in \cM_{n_\ell} \right) \eqdist \left((h_{T_{n_\ell}}(v^{(i)}) -\mu_{n_\ell})/\sigma_{n_\ell},\, i\in [k] \right).
\end{align*}
So it suffices to prove\footnote{This convergence was already implicit from the convergence of $\cP^{(n_\ell)}\todist \cP^\eps$ in \cite{AddarioEslava17}.} that $|\cM_{n_\ell}|\todist \cM_\eps$ to establish \refC{cor:max conv}.

Let $j\in \Z$. Define $X_j^{(n)}:=X_j^{(n)}(\R)\eqdist \Poi{\lambda_j}$ and $X_{\ge j+1}^{(n)}:=X_{\ge j+1}^{(n)}(\R)\eqdist \Poi{\lambda_j}$. For any $k\in \N$, as $n_\ell \to \infty$,
\begin{align*}
    \p{\Delta_{n_\ell}=\floor{\log n_\ell}+j,|\cM_{n_\ell}|=k}
    &=\p{X_j^{(n)}=k,\, X_{\ge j+1}^{(n)}=0} \\
    &\to \p{X_j=k,\, X_{\ge j+1}=0}. 
\end{align*}
It follows from the independence of $X_j^{(n)}$ and $X_{\ge j+1}^{(n)}$ that 
\begin{align}\label{eq:lim1}
    \lim_{J\to \infty}\lim_{n_\ell\to \infty} \sum_{j=-J}^J \p{\Delta_{n_\ell}=\floor{\log n_\ell}+j,|\cM_{n_\ell}|=k}
    =& \sum_{j\in \Z} e^{-2\lambda_j}\frac{\lambda_j^k}{k!}.
\end{align}
It remains to show that the tails are negligible. That is,
\begin{align}\label{eq:lim2}
    &\lim_{J\to \infty} \lim_{n_\ell\to \infty}\p{\Delta_{n_\ell}<\floor{\log n_\ell}-J}+\p{\Delta_{n_\ell}>\floor{\log n_\ell}+J}\nonumber\\
    =&\lim_{J\to \infty}\lim_{n_\ell\to \infty} \p{X_{\ge J}^{(n_\ell)}=0}+\p{X_{\ge J+1}^{(n_\ell)}>1}\nonumber \\
    =& \lim_{J\to \infty} e^{-\lambda_{-J+1}}+1-e^{-\lambda_{J}}=0.
\end{align}
Indeed, $\liminf_{n_\ell\to \infty} \p{|\cM_{n_\ell}|=k}$ is lower bounded by \eqref{eq:lim1} and, similarly, the sum of \eqref{eq:lim1} and \eqref{eq:lim2} is an upper bound for $\limsup_{n_\ell\to \infty} \p{|\cM_{n_\ell}|=k}$, completing the proof.
\end{proof}

\section*{Acknowledgements}

I would like to thank Louigi Addario-Berry and Henning Sulzbach for some very helpful discussions as well as two anonymous referees for helpful comments that substantially improved this paper.

\end{document}